\documentclass[11pt]{amsart}
\usepackage{amsmath,amsthm,amssymb,latexsym,epic,bbm,comment}
\usepackage{graphicx, stmaryrd, xcolor, tikz-cd}
\usepackage{ytableau}
\usepackage[all,2cell]{xy}
\xyoption{2cell}
\usepackage[active]{srcltx}
\usepackage[parfill]{parskip}
\usepackage[
colorlinks=true,
linkcolor=black, 
anchorcolor=black,
citecolor=black, 
urlcolor=black, 
]{hyperref}
\usepackage[e]{esvect}
\usepackage{mathrsfs}
\usepackage{todonotes}
\usepackage{enumitem}
\usepackage{booktabs}
\usepackage{makecell}
\usepackage{array}

\newtheorem{thm}{Theorem}

\newtheorem{prop}[thm]{Proposition}

\theoremstyle{definition}

\newtheorem{example}[thm]{Example}
\newtheorem{rem}[thm]{Remark}

\newtheorem{definition}[thm]{Definition}

\newcommand{\bbC}{\mathbb{C}}

\newcommand{\bbZ}{\mathbb{Z}}

\newcommand{\cC}{\mathscr{C}}

\newcommand{\cE}{\mathscr{E}}

\newcommand{\cI}{\mathscr{I}}

\newcommand{\cS}{\mathscr{S}}

\newcommand{\cA}{\mathscr{A}}

\newcommand{\rA}{\mathrm{A}}
\newcommand{\rB}{\mathrm{B}}

\newcommand{\rF}{\mathrm{F}}
\newcommand{\rG}{\mathrm{G}}
\newcommand{\rH}{\mathrm{H}}

\newcommand{\rR}{\mathrm{R}}

\newcommand{\sfA}{\mathsf{A}}

\newcommand{\sfC}{\mathsf{C}}

\newcommand{\sfv}{\mathsf{v}}

\newcommand{\A}{\mathcal{A}}

\newcommand{\C}{\mathcal{C}}

\def\H{{\mathcal{H}}}

\newcommand{\J}{\mathcal{J}}

\def\L{\mathcal{L}}

\newcommand{\X}{\mathcal{X}}
\newcommand{\Y}{\mathcal{Y}}

\newcommand{\add}{\operatorname{add}}
\newcommand{\id}{\mathsf{id}}
\newcommand{\Id}{\mathsf{Id}}

\newcommand{\bfM}{\mathbf{M}}
\newcommand{\bfK}{\mathbf{K}}

\newcommand{\bfN}{\mathbf{N}}

\newcommand{\bfC}{\mathbf{C}}

\newcommand{\tone}{\mathtt{1}}

\newcommand{\ti}{\mathtt{i}}
\newcommand{\tj}{\mathtt{j}}
\newcommand{\tk}{\mathtt{k}}
\newcommand{\tl}{\mathtt{l}}

\newcommand{\tn}{\mathtt{n}}

\newcommand{\Hom}{\mathrm{Hom}}
\newcommand{\End}{\mathrm{End}}

\newcommand{\lmod}{\text{-}\mathrm{mod}}
\newcommand{\lproj}{\text{-}\mathrm{proj}}

\newcommand{\ul}[1]{\underline{#1}}

\newcommand{\circv}{\circ_{\mathsf{v}}}
\newcommand{\circh}{\circ_{\mathsf{h}}}

\newcommand{\afmod}{\text{-}\operatorname{afmod}}
\newcommand{\sfmod}{\text{-}\operatorname{sfmod}}

\begin{document}

\title{Categorification in representation theory}

\author{Vanessa Miemietz}
\address{V.M.: School of Engineering, Mathematics and Physics, University of East Anglia, Norwich NR4 7TJ, United Kingdom,  \newline \href{https://research-portal.uea.ac.uk/en/persons/vanessa-miemietz/}{https://research-portal.uea.ac.uk/en/persons/vanessa-miemietz/}}
\email{v.miemietz@uea.ac.uk}




\begin{abstract}
In this survey article, we give an introduction to the relatively young subject of $2$-representation theory, which studies categorifications of important players in classical representation theory. In particular, we provide a streamlined exposition of the main results leading to the classification of simple $2$-representations for Soergel bimodules associated to finite Weyl groups in characteristic $0$, obtained in \cite{MMMTZ2}.
\end{abstract}

\maketitle


\section*{Introduction}

Instances of categorification have seen tremendous success in representation theory, but also other areas (in particular topology) in the last 25 years.
Categorification is generally understood as replacing objects of a certain categorical level (e.g. numbers, viewed as having "categorical level" $0$, or vector spaces) by objects of a higher categorical level (e.g. vector spaces resp. categories). 
The advantage of the higher category is that it contains more structure which can give information about the lower category. This idea  first arose in work of Crane and Frenkel \cite{CF} in mathematical physics over 30 years ago, predicting four-dimensional topological quantum field theories as categorifications of three-dimensional ones. 

One of the most famous examples is Khovanov's categorification of the Jones polynomial \cite{Kh}, which not only produced a new family  of knot and link invariants but also gave rise to a new homology theory, Khovanov homology. His conjectural link to a functorial action of the Temperley-Lieb algebra on certain categories of representations of Lie algebras \cite{BFK} was later proved by Brundan \cite{B1} and Stroppel \cite{St1}, and his categorification has been generalised in various ways, providing previously unknown knot invariants. 

At the same time, within representation theory, Lascoux, Leclerc and Thibon \cite{LLT} conjectured a relationship between decomposition numbers for Hecke algebras and canonical bases of certain Fock spaces. This was later proved by Ariki \cite{A:96} through a functorial action of an affine Lie algebra on categories of representations of affine Hecke algebras of type $A$, which also categorifies Kleshchev's modular branching rules.

Furthermore, categorifications of Hecke algebras through Soergel bimodules encode the important numerical data of Kazhdan--Lusztig polynomials \cite{Soe2}, and one of the most outstanding results in recent years was the algebraic  proof of  the Kazhdan--Lusztig conjectures by Elias and Williamson \cite{EW}, through a description of the underlying $2$-categories in the diagrammatic spirit of quiver Hecke algebras given in \cite{EW1}. Taking this approach further has led Williamson to another major breakthrough, giving counter-examples to the celebrated Lusztig and James conjectures \cite{Wi}. 

The philosophy of studying functorial actions of an algebra on a certain category is nowadays usually formulated in terms of $2$-categories.
In all cases described above, there is a $2$-category, which acts on certain categories of representations (of Lie or of Hecke algebras) and which decategorifies to a known algebra (a Temperley-Lieb, an affine Lie algebra or a Hecke algebra). The importance of the $2$-categorical structure of the functorial action of a Lie algebra on a category was first exhibited and axiomatised in the groundbreaking article on $\mathfrak{sl}_2$-categorification by Chuang and Rouquier \cite{CR}, and their study of $2$-representations of the resulting $2$-categories enabled them to prove Brou\'e's abelian defect group conjecture for blocks of the symmetric groups.

The success of $2$-representation-theoretic concepts in examples illustrates the need for an overarching theory of $2$-representations of $2$-categories, lifting the theory of representations of an algebra, as was first pointed out by  Rouquier \cite{Ro}. The development of $2$-representation theory suited to the above-mentioned examples was inititated by Mazorchuk and myself in \cite{MM1,MM2, MM3, MM5}, and later taken up by other authors. One of the major achievements in the area to date has been the classification of simple transitive $2$-representations of Soergel bimodules  in characteristic $0$ for finite Coxeter types except $H_3$ and $H_4$ in \cite{MMMTZ2}. The proof works via a reduction to certain associated `asymptotic' fusion categories defined by Lusztig, which play an important role in \cite{Lu1} and \cite{RT}. This reduction works for all Coxeter types and the gaps in the classification for $H_3$ and $H_4$ stem from a lack of understanding of some of the relevant fusion categories in those types, though significant progress has recently made in \cite{ERT}. 

In this survey article, I will provide a streamlined exposition of the ingredients needed to arrive at the main result of \cite{MMMTZ2}. In particular, I will focus on making it accessible to non-experts, by defining all the basic notions, such as $2$-categories, their $2$-representations, etc., carefully, but then only citing the more technical statements. In this way, I hope to explain the main ingredients in the proof in a way that their statements and relevance make sense, but the line of argument is not obfuscated by technical detail. In particular, I have opted to keep everything in the language of $2$-representations and $2$-representations, and completely avoided coalgebra $1$-morphisms and their (bi)comodules. In reality, the machinery of being able to realise $2$-representations as categories of comodules over coalgebra $1$-morphisms in certain completions (in the case the injective abelianisation) of the $2$-category is the most powerful tool available upon which all proofs rely. However, the way this tool is used requires frequent going back and forth across a certain biequivalence (see \cite[Theorem 4.26]{MMMTZ1}), which complicates notation in a way that was not desirable for this survey.

The article is structured as follows: In Section \ref{cats2cats}, we first explain how to view classical representations of an algebra as functors from $\Bbbk$-linear categories to $\Bbbk\lmod$ and then introduce $2$-categories, and in particular fiat $2$-categories, defining Soergel bimodules in Example \ref{soergelex}. In Section \ref{2repbasics}, we define $2$-representations and their morphisms, as well as modifications between morphisms of $2$-representations. We also introduce cell structures and cell $2$-representations and give relevant examples. In Section \ref{Hcelldoublecent}, we give the two main ingredients in the classification of simple $2$-representations for Soergel bimodules for finite Coxeter groups in characteristic zero. In Section \ref{classsec}, we introduce asymptotic algebras for Hecke algebras, their categorifications via asymptotic bicategories, and compare $2$-representations of asymptotic bicategories for a fixed $H$-cell to those of the corresponding $H$-cell reduction. We finish by collecting these ingredients into the main theorem.

\section{$\Bbbk$-linear categories and $2$-categories}\label{cats2cats}

\subsection{Classical representation theory}

Let $\Bbbk$ be an algebraically closed field.

Traditionally, an associative $\Bbbk$-algebra $A$ is a $\Bbbk$-vector space with a ring structure such that the action of scalars is central. A finite-dimensional representation  $A$ is a pair $(M,\rho)$ of a finite-dimensional vector space $M$ and an algebra homomorphism $\rho\colon A\to \End_\Bbbk(M)$. 

Alternatively, we can view an algebra over $\Bbbk$ as a $\Bbbk$-linear category $\mathsf{A}$ with finitely many objects $\mathtt{1},\dots, \mathtt{n}$. Then $A=\mathrm{End}_{\mathsf{A}}(\bigoplus_{\mathtt{i}=\mathtt{1}}^{\mathtt{n}} \mathtt{i})$ is an associative $\Bbbk$-algebra in the above sense with idempotents $e_i=\id_\ti$ such that $1_A = \sum_{1}^n e_i$. In this terminology, a representation $(M,\rho)$ of $A$ can be described as a $\Bbbk$-linear functor from $\mathsf{A}$ to the category $\Bbbk\lmod$ of finite-dimensional vector spaces, sending the object $\ti$ to the subspace $\rho(e_i)M$ and a morphism $a = e_jae_i\in \Hom_{\mathsf{A}}(\ti,\tj)$ to $\rho(a) \colon \rho(e_i)M\to\rho(e_j)M$. Conversely given a $\Bbbk$-linear functor $\mathsf{M}$ from $\mathsf{A}$ to $\Bbbk\lmod$, we obtain a representation in the above classical sense by setting $M = \bigoplus_{\ti = \tone}^\tn \mathsf{M}(\ti)$ defining $\rho(a) = \mathsf{M}(a)$.

This explains why, in order to categorify representations of algebras, we will need to consider $2$-categories and $2$-functors to a fixed target category.

Notice that since $\Bbbk\lmod$ is an additive category, the category of $\Bbbk$-linear functors from $\mathsf{A}$ to $\Bbbk\lmod$ is naturally equivalent to the category of $\Bbbk$-linear functors from the additive envelope of $\mathsf{A}$ to $\Bbbk\lmod$, where the additive envelope of $\mathsf{A}$ denotes the category whose objects are tuples of objects in $\mathsf{A}$ and morphisms are given by matrices of morphisms in $\mathsf{A}$.

We say a category $\sfC$ is {\bf finitary} (over $\Bbbk$) if it is small, $\Bbbk$-linear, additive, idempotent complete (also known as Karoubi complete or Cauchy complete), has only finitely many indecomposable objects up to isomorphism and finite-dimensional morphism spaces. In other words, $\sfC$ is finitary if and only if it is equivalent to the category $A\lproj$ of finitely generated projective $A$-modules for a finite-dimensional associative $\Bbbk$-algebra $A$.

\subsection{$2$-categories}

We denote by $\boldsymbol{1}$ the trivial category consisting a single object $\star$ and its identity morphism $\id_\star$.

\begin{definition}
A {\bf $2$-category $\cC$} is a category enriched over the monoidal category $\mathbf{Cat}$ of small categories, i.e. it consists of 
\begin{itemize}
 \item a class $\cC$ of objects;
\item  for every $\mathtt{i},\mathtt{j}\in\cC$ a small category $\cC(\mathtt{i},\mathtt{j})$ of morphisms from
$\mathtt{i}$ to $\mathtt{j}$ where
\begin{itemize}
\item objects in $\cC(\mathtt{i},\mathtt{j})$ are called {\bf $1$-morphisms},
\item morphisms in $\cC(\mathtt{i},\mathtt{j})$ are called {\bf $2$-morphisms};
\end{itemize}
\item functorial composition $-\circ- \colon \cC(\mathtt{j},\mathtt{k})\times \cC(\mathtt{i},\mathtt{j})\to
\cC(\mathtt{i},\mathtt{k})$;
\item a functor $\mathbbm{1}_\ti \colon \boldsymbol{1} \to \cC(\ti,\ti)$;
\end{itemize}
such that
\begin{itemize}
\item the diagram of functors 
\[ \xymatrix{
\cC(\mathtt{k},\mathtt{l})\times\cC(\mathtt{j},\mathtt{k})\times \cC(\mathtt{i},\mathtt{j}) \ar^{(-\circ-)\times \Id}[rr] \ar_{\Id\times (-\circ-)}[d]&&\cC(\mathtt{j},\mathtt{l})\times \cC(\mathtt{i},\mathtt{j})\ar^{-\circ-}[d]\\
\cC(\mathtt{k},\mathtt{l})\times\cC(\mathtt{i},\mathtt{k})\ar^{-\circ-}[rr]&& \cC(\mathtt{i},\mathtt{l})
}\]
commutes;
\item the diagrams of functors
\[ \xymatrix{
\boldsymbol{1} \times \cC(\ti,\tj) \ar_{\sim}[dr]\ar^{\mathbbm{1}_\tj \times \Id}[r] &\cC(\tj,\tj)\times \cC(\ti,\tj)\ar^{-\circ-}[d]\\
& \cC(\ti,\tj)
} \quad\text{and}\quad  \xymatrix{
 \cC(\ti,\tj) \times \boldsymbol{1}\ar_{\sim}[dr]\ar^{\Id\times \mathbbm{1}_\ti }[r] &\cC(\ti,\tj)\times \cC(\ti,\ti)\ar^{-\circ-}[d]\\
& \cC(\ti,\tj)
}\]

commute.
\end{itemize}
\end{definition}
By slight abuse of notation, we will denote the $1$-morphism $\mathbbm{1}_\ti(\star)\in \cC(\ti,\ti)$ simply by $\mathbbm{1}_\ti$.

Writing $-\circh -$ for the so-called {\bf horizontal} composition of morphisms under the functor $-\circ -$ and $-\circv-$ for the so-called {\bf vertical} composition of $2$-morphisms inside each morphism category, functoriality of $-\circ -$ in particular implies that for $1$-morphisms $\rF, \rF',\rF''\in \cC(\ti,\tj)$ and $\rG,\rG',\rG''\in \cC(\tj,\tk)$, as well as $2$-morphisms $\alpha \colon \rF\to\rF', \alpha' \colon \rF'\to\rF'', \beta \colon \rG\to\rG', \beta' \colon \rG'\to\rG''$, we have 
\[ (\beta'\circv \beta)\circh(\alpha'\circv \alpha) = (\beta'\circh\alpha' )\circv(\beta\circh \alpha), \]
the so-called interchange law.

\begin{rem}
If we only require the diagrams of functors in the definition of a $2$-category to commute up to coherent natural isomorphisms, we obtain the notion of a bicategory. We will encounter some bicategories later in this article. However, since the bicategorical structure of these is not going to be used explicitly, we will ignore the technicalities for a more accessible exposition. Full details and proofs for bicategories of the relevant statements can be found in \cite{MMMTZ1}.
\end{rem}

\begin{example}
\begin{enumerate}[label=\alph*)]
\item The $2$-category $\mathbf{Cat}$ has
\begin{itemize}
\item small categories as objects;
\item functors as $1$-morphisms ;
\item natural transformations as $2$-morphisms.
\end{itemize}
\item The $2$-category $\mathfrak{A}_{\Bbbk}^{f}$ has
\begin{itemize}
\item finitary categories as objects;
\item $\Bbbk$-linear (thus automatically additive) functors as $1$-morphisms;
\item all natural transformations of such functors $2$-morphisms.
\end{itemize}
\end{enumerate}
\end{example}

The relevance of the second example stems from the fact that this is what the examples of categorifications of finite-dimensional algebras in representation theory are enriched over. In particular, this implies that it naturally plays the role of $\Bbbk\lmod$ in classical representation theory, and should be the target of any sensible $2$-representation.

\begin{definition}
A $2$-category $\cC$ is {\bf finitary} over $\Bbbk$ if
\begin{itemize}
\item $\cC$ has finitely many objects;
\item each $\cC(\mathtt{i},\mathtt{j})$ is  in $\mathfrak{A}_{\Bbbk}^{f}$;
\item composition is {$\Bbbk$-bilinear} (and hence biadditive);
\item identity $1$-morphisms are {indecomposable}.
\end{itemize}
\end{definition}

Finitary $2$-categories should be thought of as $2$-analogues of finite-dimensional algebras. Indeed, $\cC$ having only finitely many objects and each $\cC(\ti,\tj)$ having only finitely many indecomposables up to isomorphism guarantees that the decategorification of $\cC$, i.e.\ its split Grothendieck group, is a finite-dimensional algebra (over $\mathbb{Z}$). Indecomposability of $1$-morphisms is not crucial, in fact, it is easy to go back and forth between a finitary $2$-category and its additive closure, which only has one object, see \cite[Section~2.4]{MMMTZ1}, but it makes life easier in the same way that knowing a decomposition of the identity of an algebra into a sum of primitive orthogonal idempotents does.

\begin{definition}
A $2$-category $\cC$ is {\bf fiat}  ({\bf f}initary~- {\bf i}nvolution~- 
{\bf a}djunction~- {\bf t}wo-category) if
\begin{itemize}
\item it is finitary; 
\item there is a weak involutive equivalence $(-)^*\colon\cC \to \cC^{\mathrm{co},\mathrm{op}}$ (i.e.\ the $2$-category, in which the directions of both $1$-morphisms and $2$-morphisms are reversed) such that, for any objects $\ti,\tj\in \cC$ and any $\rF\in \cC(\ti,\tj)$, there exist
adjunction morphisms $\rF\circ \rF^*\rightarrow \mathbbm{1}_{\mathtt{j}}$ and
$\mathbbm{1}_{\mathtt{i}}\rightarrow \rF^*\circ \rF$.
\end{itemize}
\end{definition}

\begin{example}\label{CAex} Let $A$ be a connected finite-dimensional $\Bbbk$-algebra. The $2$-category $\cC_A$ has 
\begin{itemize}
\item one object $\bullet$ (identified with $A$-$\mathrm{proj}$);
\item $1$-morphisms are endofunctors of $\varnothing$ isomorphic to tensoring with bimodules in the additive closure of $A\oplus A\otimes_\Bbbk A$;
\item $2$-morphisms are natural transformations of such functors (corresponding to bimodule momorphisms).
\end{itemize}

We observe that $\cC_A$ is always finitary. If $A$ is basic  with complete set of idempotents $e_1,\dots, e_n$, the indecomposable $1$-morphisms correspond to the bimodules $A$ and $Ae_i\otimes_\Bbbk e_jA$, for $i,j=1,\dots n$. If $A$ is weakly symmetric, $\cC_A$ is fiat with involution given by $(Ae_i\otimes_\Bbbk e_jA)^*\cong Ae_j\otimes_\Bbbk e_iA$.
\end{example}

\begin{example}\label{soergelex} \emph{(Soergel bimodules)} Let $(W,S)$ be a Coxeter group, i.e. $S$ is a finite set and $W$ is the group given by
\[W=\langle s_i \vert \, s_i\in S, s_i^2=1, (s_is_j)^{m_{ij}}=1\rangle,\]
for some $m_{ij}\in \mathbb{Z}_{\geq 2}\cup \{\infty\}$. Let $V$ be a reflection faithful representation of $W$. Define the coinvariant algebra $\rR=\bbC[V]/(\bbC[V]^W)_+$ where $\bbC[V]$ is the symmetric algebra of $V$ and $(\bbC[V]^W)_+$ is the ideal generated by the invariants of positive degree under the natural action of $W$ on $\bbC[V]$. We view $\rR$ as a graded algebra where the degree of (the coset of) a homogeneous polynomial is twice its polynomial degree.
Setting $\rR_i:=\rR\otimes_{\rR^{s_i}}\rR$, the category of {\bf Soergel bimodules} is the additive closure of finite tensor products (over $R$) of grading shifts of the $R_i$.

We define the $2$-category $\cS=\cS_{W,S,V}$ of to have 
\begin{itemize}
\item one object $\varnothing$, which we identify with the category of finitely generated graded projective $\rR$-modules  $\rR$-$\mathrm{gproj}$;
\item as $1$-morphisms all endofunctors of $\varnothing$ which are isomorphic to tensoring with Soergel bimodules;
\item as $2$-morphisms all graded natural transformations of such functors (which correspond to bimodule morphisms of degree zero). 
\end{itemize}

If $W$ is finite, then $\cS$ is $\bbZ$-fiat (meaning it satisfies all conditions of the definition of fiat, except for having infinitely many isomorphism classes of indecomposable $1$-morphisms, but only finitely many up to grading shift), and all tools of finitary representation theory apply (c.f.\ \cite[Subsections 2.7, 3.1]{MMMTZ2} - note that in loc.cit. our $\cS$ was called $\cS^{(0)}$). Moreover, $\cS$ categorifies the Hecke algebra $\mathsf{H} = \mathsf{H}_{W,S}$ by \cite{Soe2}, see also \cite{EW}. 

An excellent reference for all things specific to the example of Soergel bimodules in this article, e.g.\ Coxeter groups, reflection faithful representations, Hecke algebras and Soergel bimodules themselves is the book by Elias--Makisumi--Thiel--Williamson \cite{EMTW}. However, we remark that therein Soergel bimodules are considered over $\bbC[V]$ rather than over the coinvariant algebra, as in the original definition by Soergel. This results in the corresponding category of Soergel bimodules not being finitary. Since the total invariants $(\bbC[V]^W)_+$ generate a central $2$-ideal that necessarily acts by zero in any simple $2$-representation, this does not affect the $2$-representation theoretic questions considered here and the classification of simple $2$-representations using Soergel's definition and the one used in \cite{EW} and \cite{EMTW} is the same.
\end{example} 

In analogy to needing ideals to take quotients of algebras, we need $2$-ideals to be able to take quotients of $2$-categories

\begin{definition}
A $2$-ideal $\cI$ in a $2$-category $\cC$ consists of, for each $\ti,\tj\in \cC$ and each $\rF,\rG\in \cC(\ti,\tj)$, a subspace $\Hom_{\cI(\ti,\tj)}(\rF,\rG)\subseteq\Hom_{\cC(\ti,\tj)}(\rF,\rG)$ such that the collection of these subspaces is closed under both horizontal and vertical composition with any $2$-morphism in $\cC$.
\end{definition}

\section{Basic definitions of $2$-representation theory}\label{2repbasics}
\subsection{$2$-representations}

\begin{definition}
A {\bf finitary $2$-representation} $\bfM$ of a finitary $2$-category $\cC$ is a $2$-functor $\cC \to \mathfrak{A}_{\Bbbk}^f$, i.e.
\begin{itemize}
\item for any $\ti \in \cC$, we have $\bfM(\ti)\in  \mathfrak{A}_{\Bbbk}^f$;
\item for any $\ti,\tj\in\cC$, there $\bfM$ defines a $\Bbbk$-linear functor $\cC(\ti,\tj)$ to the category of $\Bbbk$-linear functors from $\bfM(\ti)$ to $\bfM(\tj)$;
\item for any $\ti \in \cC$, we have $\bfM(\mathbbm{1}_\ti) = \Id_{\bfM(\ti)}$;
\item for any $\rF\in \cC(\ti,\tj), \rG\in \cC(\tj,\tk)$, we have $\bfM(\rG\rF)=\bfM(\rG)\bfM(\rF)$.
\end{itemize}
\end{definition}
In particular, this means that or any $\ti \in \cC$, the category $\bfM(\ti)$ is equivalent to $B_{\ti}$-$\mathrm{proj}$ for some algebra $B_{\ti}$, for any $1$-morphism $\rF\in \cC(\ti,\tj)$, there is a $\Bbbk$-linear functor $\bfM(\rF)\colon \bfM(\ti)\to\bfM(\tj)$  and, for any $2$-morphism $\alpha\colon \rF\to\rG$, there  is a natural transformation $\bfM(\alpha)\colon \bfM(\rF)\to\bfM(\rG)$. 
Finitary $2$-representations  of $\cC$ together with $2$-natural transformations and 
modifications form a $2$-category, denoted by $\cC\afmod$, see e.g.\ \cite[Proposition~1]{MM3}. 

\begin{example}
\begin{enumerate}[label=\alph*)]
\item For each $\mathtt{i}\in\cC$, we have the so-called {\bf principal}  $2$-representation $\mathbf{P}_\mathtt{i}= \cC(\mathtt{i}, -)$. 
\item $\cC_A$ and $\cS$ were defined via their {\bf  natural} $2$-representations on $A$-$\mathrm{proj}$, resp.\ $R$-$\mathrm{proj}$.
\end{enumerate}
\end{example}


An appropriate analogue of a simple representation of an algebra was defined in \cite{MM5}. In loc.\ cit.\ this was called a simple transitive $2$-representation, but we will opt for the terminology of just calling it a simple $2$-representation, since transitivity is implied by simplicity. 

\begin{definition}
A finitary $2$-representation $\bfM$ of a finitary $2$-category $\cC$ is called {\bf simple} if 
$\coprod_{\ti\in\cC}\mathbf{M}(\ti)$ has no proper $\cC$-stable ideals.
\end{definition}

The terminology is justified by the fact that this notion of simple $2$-representation admits a weak version of the classical Jordan--H\"older theorem, see \cite[Theorem 8]{MM5}.

\begin{rem}
In the same way that the definition of $2$-category can be relaxed to that of a bicategory, the notion of $2$-representation can be relaxed to obtain a notion of birepresentation. This is then a $\Bbbk$-linear pseudofunctor from the $2$-category (or bicategory) one wishes to represent to the chosen target. Fixing the same target $ \mathfrak{A}_{\Bbbk}^f$, we can thus speak about finitary or simple birepresentations. Using a result by Powers \cite[Section 4.2]{Pow}, any finitary birepresentation of a finitary $2$-category $\cC$ is equivalent to a $2$-representation, so we will only use the notion of birepresentation when considering bicategories that are not $2$-categories. Again, we refer the reader to \cite{MMMTZ1} for details about the bicategorical setting.
\end{rem}

\begin{rem}
Any monoidal category can be viewed as a $2$-category with one object, whose endomorphism category is given by the monoidal category in question. Under this identification, a $2$-representation corresponds to a module category in the language of \cite{EGNO}. In the special case where the monoidal category is a fusion category, the simple $2$-representations are precisely the indecomposable module categories studied in \cite{EGNO}.
\end{rem}

\subsection{Injective abelianisations}

Let $\C$ be a finitary category. We denote by $\ul{\C}$ its {\bf injective abelianisation}, which has diagrams 
$X_0\xrightarrow{f}X_1$ as objects and whose morphisms are given by pairs $(\phi_0,\phi_1)$ such that the solid part of the diagram
\[\xymatrix{
X_0\ar^f[rr]\ar^{\phi_0}[d]&&X_1\ar^{\phi_1}[d] \ar_{\eta}@{.>}[dll]\\
Y_0\ar^g[rr]&&Y_1\\
}\]
commutes, modulo the $\Bbbk$-subspace spanned by those pairs $(\phi_0,\phi_1)$ such that there exists an $\eta$ as indicated by the dotted arrow with $\phi_0 = \eta f$. This is an abelian category whose injective objects are precisely the diagrams $X\to 0$, for $X\in\C$, see \cite{Fr}.

We can abelianise 
a finitary $2$-representation $\bfM$ of a finitary $2$-category $\cC$ by doing so locally, i.e. setting 
$\ul{\bfM}(\ti) = \ul{\bfM(\ti)}$ and defining the action of $\cC$ component-wise.

\subsection{Cells and cell $2$-representations}

One important observation in \cite{MM1} was that $2$-categories come with a natural cell theory that generalises the cell theory classically used to study semigroups and monoids, see \cite{Gr}.

On (isomorphism classes of) indecomposable $1$-morphisms in a finitary $2$-category $\cC$, define 
the
{\bf left preorder} by saying $\rF\leq_L \rG$ if there exists an $\rH$ such that $\rG$ is isomorphic to a direct summand of $\rH\rF$. Similarly we define the {\bf right preorder}  by saying $\rF\leq_R \rG$ if there exists an $\rH$ such that $\rG$ is isomorphic to a direct summand of $\rF\rH$, as well as the {\bf two-sided preorder} by saying that $\rF\leq_J \rG$ if there exist  $\rH_1,\rH_2$ such that $\rG$ is isomorphic to a direct summand of $ \rH_1\rF\rH_2$. Equivalence classes with respect to the left, respectively right or two-sided, preorders are called {\bf left}, respectively {\bf right} or {\bf two-sided}, {\bf cells}.

An {\bf $H$-cell} is the intersection of a left and a right cell. If $\cC$ is fiat, then for any left cell $\L$, we have the dual right cell $\L^*=\{\rF^*\,\mid\, \rF\in \L\}$. The intersection $\L\cap \L^*$ is called a {\bf diagonal $H$-cell}.

Fix a left cell  $\L$ in $\cC$ and note that, by construction, there is a unique object $\ti_{\L}\in \cC$, such that any $\rF\in L$ belongs to $\cC(\ti_{\L}, \tj)$ for some $\tj\in \cC$. For the diagonal $H$-cell $\H_{\L} = \L\cap\L^*$, this implies that every $\rF\in \H_{\L}$ belongs to $\cC(\ti_{\L},\ti_{\L})$. Pick $\rF\in \H_{\L}$ and consider the $2$-representation $\bfN_{\L}$ given by $\bfN_{\L}(\tj) = \add_{\cC(\ti,\tj)}\{\rH\rF\,\mid\, \rH\in \cC(\ti_{\L},\tj)\}$ with the natural $\cC$-action by multiplication on the left. This has a unique maximal $\cC$-stable ideal (see e.g.\ \cite[Subsection 6.5]{MM2}), by which we can quotient to obtain a simple $2$-representation, the so-called  {\bf cell $2$-representation} $\mathbf{C}_{\L}$.

By \cite[Subsection 3.2]{ChMa}, any simple $2$-representation $\bfM$ of a finitary $2$-category $\cC$ has a so-called {\bf apex},
which is the unique maximal two-sided cell $\J$ of $\cC$ with respect to the property of not being annihilated by $\bfM$. Thus, if we want to classify simple $2$-representations for a given $\cC$, we can do so apex by apex. By construction, for a left cell $\L$ in a two-sided cell $\J$, the apex of $\mathbf{C}_{\L}$ is $\J$.

\begin{example}\label{cellex}
\begin{enumerate}[label=\alph*)]
\item
 For $\cC_A$ (c.f.\ Example \ref{CAex}), the cells are given by

\begin{center}

\begin{tabular}{|c|}
\hline
$\mathbbm{1}=A$\\
\hline
\end{tabular}

\vspace{3mm}
{\renewcommand{\arraystretch}{2.5}

\begin{tabular}{|c|c|c|c|}
\hline
$Ae_1\otimes e_1A$ & $Ae_1\otimes e_2A$ & $\cdots$ & $Ae_1\otimes e_nA$\\
\hline
$Ae_2\otimes e_1A$ & $Ae_2\otimes e_2A$ & $\cdots$ & $Ae_2\otimes e_nA$\\
\hline
$\vdots$ & $\vdots $& $\ddots$ & $\vdots$\\
\hline
$Ae_n\otimes e_1A$ & $Ae_n\otimes e_2A$ & $\cdots$ & $Ae_n\otimes e_nA$\\
\hline
\end{tabular}
}
\end{center}
where the big squares are two-sided cells, the columns are left cells and the rows are right cells. In this example, each $H$-cell only contains one element, a property that we call {\bf strongly regular} and which has strong $2$-representation theoretic implications, such as the cell $2$-reprsentations for all left cells in a two-sided $\J$ being equivalent, and being the only simple $2$-representations with apex $\J$ up to equivalence.  

\item\label{Scellex} For $\cS$ (c.f.\ Example \ref{soergelex}), indecomposable $1$-morphisms in $\cS$ are labelled by  $W$ and denoted by $\theta_w, w\in W$. More precisely, if $w= s_{i_1}\cdots s_{i_m}$ is a reduced expression, $\theta_w$ is the unique indecomposable summand of $R_{i_1}\otimes_R\cdots \otimes_RR_{i_m}$ which does not appear as a direct summand in any $R_{j_1}\otimes_R\cdots \otimes_RR_{j_l}$ for $l<m$.
Appropriate shifts of the $\theta_w$ descend to the Kazhdan--Lusztig basis under decategorification and the cell structure of $\cS$ corresponds to Kazhdan--Lusztig cells. If $W$ is the dihedral group with $8$ elements, i.e. 
\[W=\langle s,t \vert s^2=1=t^2, stst=tsts \rangle,\]
cells are of the form
\begin{center}

\begin{tabular}{|c|}
\hline
$\mathbbm{1}$\\
\hline
\end{tabular}

\vspace{3mm}

\begin{tabular}{|c|c|}
\hline
$\theta_s, \theta_{sts}$ & $\theta_{st}$\\
\hline
$\theta_{ts}$ & $\theta_t,\theta_{tst}$\\
\hline
\end{tabular}

\vspace{3mm}

\begin{tabular}{|c|}
\hline
$\theta_{stst}$\\
\hline
\end{tabular}
\end{center}
We remark that the diagonal $H$-cells in the middle two-sided cell contain two elements. By slight abuse of terminology, we will say $v$ and $w$ are related or equivalent with respect to one of the orders if $\theta_v$ and $\theta_w$ are. This terminology coincides with the original terminology for Kazhdan--Lusztig orders and cells.
\end{enumerate}
\end{example}

\subsection{Morphisms of $2$-representations and modifications}

Throughout this subsection, let $\cC$ be a finitary $2$-category and let $\bfM$ and $\bfN$ be  $2$-representations, i.e.\ $2$-functors from $\cC$ to $\mathfrak{A}_{\Bbbk}^f$. 

\begin{definition}
A {\bf morphism $\Phi\colon \bfM\to \bfN$ of $2$-representations} is a strong natural transformation between $2$-functors. Concretely, this means That $\Phi$ consists of a collection of $\Bbbk$-linear functors ($1$-morphisms in $\mathfrak{A}_{\Bbbk}^f$) $$\Phi_\ti\colon \bfM(\ti)\to\bfN(\ti),$$ for $\ti\in \cC$, as well as a collection of natural isomorphisms ($2$-morphisms in $\mathfrak{A}_{\Bbbk}^f$) $$\Phi_\rF\colon \bfN(\rF)\Phi_\ti\to \Phi_\tj\bfM(\rF),$$ for each $\rF\in \cC(\ti,\tj)$, which satisfy what is called lax naturality, namely that for each $\rF\in \cC(\ti,\tj),\rG\in \cC(\tj,\tk)$, we have 
\[\xymatrix{
\bfM(\ti)\ar^-{\bfM(\rG\rF)}[rr]\ar_-{\Phi_\ti}[dd]&&\bfM(\tk) \ar^-{\Phi_\tk}[dd]&&\bfM(\ti)\ar^-{\bfM(\rF)}[r]\ar_-{\Phi_\ti}[dd]&\bfM(\tj)\ar^-{\bfM(\rG)}[r]\ar_-{\Phi_\tj}[dd]&\bfM(\tk) \ar^-{\Phi_\tk}[dd]\\
&&&=&&&\\
\bfN(\ti)\ar^-{\bfN(\rG\rF)}[rr] 
\ar@{=>}^-{\Phi_{\rG\rF}}[uurr]
&&\bfN(\tk) &&  \bfN(\ti)\ar^-{\bfN(\rF)}[r]\ar@{=>}^-{\Phi_{\rF}}[uur]&\bfN(\tj)\ar^-{\bfN(\rG)}[r] \ar@{=>}^-{\Phi_{\rG}}[uur]&\bfN(\tk).} 
\]
\end{definition}

We can compose morphisms $\Phi\colon \bfM\to\bfN$ and $\Psi\colon \bfN\to \bfK$ of finitary $2$-representations by setting $(\Psi\Phi)_\ti = \Psi_\ti\Phi_\ti $ and $(\Psi\Phi)_\rF = (\id_{\Psi_\tj}\circh\Phi_\rF)\circv(\Psi_\rF\circh \id_{\Phi_\ti})$. That is, in functor notation, $(\Psi\Phi)_\rF \colon \bfK(\rF)\Psi_\ti\Phi_\ti\to \Psi_\tj\Phi_\tj\bfM(\rF)$ is the natural transformation whose component on an object $X\in \bfM(\ti)$ is given by the composition
\[  \bfK(\rF)\Psi_\ti\Phi_\ti (X)  \xrightarrow{(\Psi_\rF)_{\Phi_\ti(X)}}  \Psi_\tj\bfN(\rF)\Phi_\ti (X)  \xrightarrow{\Psi_\tj((\Phi_\rF)_X)} \Psi_\tj\Phi_\tj\bfM(\rF)(X) .\]

\begin{definition}
Given morphisms $\Phi, \Phi' \colon \bfM\to\bfN$ of finitary $2$-representations, a {\bf modification} $\varphi\colon \Phi\to\Phi'$ consists of the data of a $2$-morphism (i.e. a natural transformations of functors) $\varphi_\ti\colon \Phi_\ti \to \Phi'_\ti$ for every object $\ti\in \cC$ such that, for every $\rF\in \cC(\ti,\tj)$, the diagram 
\[\xymatrix{
\bfN(\rF)\Phi_\ti \ar^-{\id_{\bfN(\rF)} \circh \varphi_\ti}[rr]\ar_-{\Phi_\rF}[d]&& \bfN(\rF)\Phi'_\ti\ar^-{\Phi'_\rF}[d]\\
\Phi_\tj \bfM(\rF)\ar^-{\varphi_\tj\circh \id_{\bfM(\rF)}}[rr]&& \Phi'_\tj \bfM(\rF)
}\]
commutes.
\end{definition}

In this way, we obtain a $2$-category $\cC\afmod$ whose 
\begin{itemize}
\item objects are finitary $2$-representations of $\cC$;
\item $1$-morphisms are morphisms of finitary $2$-representations;
\item $2$-morphisms are modifications.
\end{itemize}
Moreover, we can consider, for any fixed $2$-representation $\bfM$ of $\cC$, the $2$-category $\cE nd_{\cC}(\bfM)$ whose single object is indentified with $\bfM$, with $1$-morphisms given by endomorphisms of $\bfM$ and $2$-morphisms given by all modifications between these.


Any morphism $\Phi\colon \bfM \to \bfN$ extends to the injective abelianisations by component-wise application, and we denote the resulting morphism $\underline{ \bfM} \to \underline{\bfN}$ again by $\Phi$.

\begin{definition} We call $\Phi$ {\bf exact} if its extension to the abelianisation is locally given by exact functors, and we say $\Phi$ is {\bf injective} if the extension of each $\Phi$ to the injective abelianisation is injective in the categories of functors from $ \underline{ \bfM}(\ti)$ to $\underline{\bfN}(\ti)$.
\end{definition}

\section{$H$-cell reduction and the double centraliser theorem}\label{Hcelldoublecent}
\subsection{$H$-cell reduction}\label{Hcellredsec}

One of the most powerful tools in classifying simple $2$-representations of a given $2$-category is the so-called $H$-cell reduction, developed in \cite{MMMZ}.

To define this, recall that a $2$-ideal of a $2$-category $\cC$ consists of a collection of ideals $\cI(\ti,\tj)$ in each $\cC(\ti,\tj)$, such that 
for any 
$\alpha\in \cC(\tk, \tl), \beta \in  \cI(\tj, \tk), \gamma \in \cC(\ti, \tj)$, the horizontal composition satisfies  $\alpha\circh\beta\circh\gamma \in \cI(\ti,\tl)$.

Let  $\H$ be a diagonal $H$-cell in a fiat $2$-category $\cC$, contained in a two-sided cell $\J$ and denote by $\ti = \ti_\H$ the source and target of $1$-morphisms in $\H$. We construct a new fiat $2$-category $\cC_\H$ in several steps:

\begin{itemize}
\item We first take the quotient of $\cC$ by the $2$-ideal generated by the identities on all $1$-morphisms belonging to two-sided cells $\J' \nleq \J$.
\item Inside this quotient, we let $\cC_{(\H)}$ be the $2$-full sub-$2$-category with single object $\ti$ and $1$-morphisms in the additive closure of $\mathbbm{1}_{\ti}$ and all $1$-morphisms $\rF$  in $\H$.
\item Finally, we let $\cC_\H$ be the quotient of $\cC_{(\H)}$ by the maximal $2$-ideal not containing $\id_{\rF}$ for any $\rF\in \H$.
\end{itemize}

\begin{example}\label{Hcellex}
\begin{enumerate}[label=\alph*)]
\item For $\cC=\cC_A$, take $\H = \{Ae_1\otimes e_1A\}$, then $\cC_\H$ has cell structure 
\begin{center}

\begin{tabular}{|c|}
\hline
$\mathbbm{1}=A$\\
\hline
\end{tabular}

\vspace{3mm}

\begin{tabular}{|c|}
\hline
$Ae_1\otimes e_1A$ \\ \hline
\end{tabular}
\end{center}

\item\label{SHcellex} For  $\cS=\cS_{B_2}$, take $\H = \{\theta_s, \theta_{sts}\}$, then $\cS_\H$ has cell structure

\begin{center}

\begin{tabular}{|c|}
\hline
$\mathbbm{1}=\theta_1$\\
\hline
\end{tabular}

\vspace{3mm}

\begin{tabular}{|c|}
\hline
$\theta_s, \theta_{sts}$ \\ \hline
\end{tabular}
\end{center}

\end{enumerate}
\end{example}

Denote by $\cC\sfmod_\J$ the $1$-full and $2$-full sub-$2$-category of $\cC\afmod_\J$ consisting of simple $2$-representations with apex $\J$.

\begin{thm}[{\cite[Theorem 15]{MMMZ}\cite[Theorem 4.32]{MMMTZ1}}]
Let $\H$ be a diagonal $H$-cell contained in a $2$-sided cell $\J$ of a fiat $2$-category $\cC$. There is a bijection between simple $2$-representations of $\cC$ with apex $\J$ and simple $2$-representations of $\cC_\H$ with apex $\H$.

Furthermore, this bijection is part of a biequivalence 
$$\cC\sfmod_\J \simeq \cC_\H\sfmod_\H.$$
\end{thm}

Recalling that, in order to classify simple $2$-representations, it suffices to do so apex by apex, the theorem now implies that we can pick a diagonal $H$-cell in each $2$-sided cell and classify simple $2$-representations for each $\cC_\H$. The obvious advantage of $\cC_\H$ is that it is significantly smaller, having at most one $H$-cell which does not contain an identity.

\subsection{Double Centraliser Theorem}

For any finitary $2$-representation $\bfM$ of a finitary $2$-category $\cC$, we denote by $\cE nd_{\cC}(\bfM)$ its endomorphism $2$-category, which has $\bfM$ as its sole object, whose $1$-morphisms are all strong transformations from $\bfM$ to itself and whose $2$-morphisms are all modifications between such strong transformations. 

We denote by $\cE nd^{\mathrm{inj}}_{\cC}(\bfM)$ the $2$-full sub-$2$-category whose $1$-morphisms are only those strong transfomations $\Psi$
 whose component functors $\Psi_\ti\colon \bfM(\ti)\to \bfM(\ti)$, when extended to $ \ul{\bfM}(\ti)$, are injective in the category of endofunctors of the latter.

Observe that, by construction, there is a canonical $2$-functor
$$\mathrm{can}\colon \cC_\H  \to \cE nd_{\cE nd_{\cC_\H}(\bfM)}(\bfM).$$
%
%
\begin{thm}[{\cite[Theorem 5.3]{MMMTZ1}}]
Let $\cC$ be a fiat $2$-category and $\H$ a diagonal $H$-cell in $\cC$.
Let $\bfM$ be a non-trivial simple $2$-representation of $\cC_\H$. 

The canonical $2$-functor restricts to an equivalence  $$\mathrm{add}(\H) \simeq  \cE nd^{\mathrm{inj}}_{\cE nd_{\cC_\H}(\bfM)}(\bfM).$$
\end{thm}

In good circumstances, such as those we will encounter in the next section, we obtain a bijection between simple $2$-representations of $\cC_\H$ with apex $\H$ and simple $2$-representations of $\cE nd_{\cC_\H}(\bfM)$.

\section{Classification of simple $2$-representations for $\cS$ in finite Weyl type}\label{classsec}

Throughout Section \ref{classsec}, let $(W,S)$ be a finite Coxeter group and assume that $\Bbbk=\bbC$. 

\subsection{Hecke algebras and asymptotic Hecke algebras}

The {\bf Hecke algebra} $\mathsf{H} = \mathsf{H}_{W,S}$ is the $\bbZ[\sfv,\sfv^{-1}]$-algebra generated by the Kazhdan--Lusztig generators $b_{s_i}, s_i\in S$, subject to relations 
\[ b_{s_i}^2 =  (\sfv+\sfv^{-1})b_{s_i} \qquad \text{and} \qquad \underbrace{b_{s_i}b_{s_j}b_{s_i} \cdots }_{{m_{ij}}\text{ terms}}= \underbrace{b_{s_j}b_{s_i}b_{s_j} \cdots }_{{m_{ij}}\text{ terms}}.\]
This has a Kazhdan--Lustig basis $\{ b_w \,\mid\, w\in W\}$ and the structure constants with respect to this basis are given by the Kazhdan--Lusztig polynomials
\[b_vb_w = \sum_{u\in W} h_{v,w,u}b_u\]
where $h_{v,w,u}\in \bbZ_{\geq 0}[\sfv,\sfv^{-1}]$. 
The Kazhdan--Lusztig left and right partial preorders on $W$ are defined by
\begin{align*}
v\leq_L w  &\quad\text{ if }\exists u\in W \text{ with } h_{u,v,w}\neq 0\\
v\leq_R w  &\quad\text{ if }\exists u\in W \text{ with } h_{v,u,w}\neq 0
\end{align*}
and the two-sided partial preorder is their union. The resulting cells are called left, right, resp.\ two-sided Kazhdan--Lusztig cells. 
There exists a unique $\mathbf{a}(u)$ such that $h_{v,w,u}\in \sfv^{\mathbf{a}(u)}\bbZ_{\geq 0}[\sfv^{-1}]$ for all $v,w\in W$, and the function $u\mapsto \mathbf{a}(u)$ is called Lusztig's $\mathbf{a}$-function.
The categorification theorem \cite{Soe2, EW}, see also \cite[Theorem 5.24 ]{EMTW}, implies that $\rB_w:= \theta_w\langle\mathbf{a}(w) \rangle$ decategorifies to $b_w$ and hence the coefficient of $\sfv^j$ in $h_{v,w,u}$ coincides with the multiplicity of $\rB_u\langle -j \rangle$ in $\rB_v\rB_w$ (c.f. Example \ref{cellex}\ref{Scellex}) and that the cell structure on $1$-morphisms of $\cS$ indeed coincides with Kazhdan--Lusztig's cell structure on Hecke algebras.

Taking $\gamma_{v,w,u^{-1}}$ to be the coefficient of $\sfv^{\mathbf{a}(u)}$ in $h_{v,w,u}$, we can use these to define a new associative $\bbZ$-algebra $\mathsf{A}$ on basis $\mathsf{t}_w, w\in W$, called the {\bf asymptotic Hecke algebra}. There is an algebra isomorphism 
\[\sfA = \prod_\J \sfA_\J\]
where $\J$ runs over all two-sided Kazhdan--Lusztig cells (which coincide with the two-sided cells of $\cS$) and $A_\J$ is the $\bbZ$-span of all $w\in \J$ .
For any diagonal $H$-cell $\H$, contained in a two-sided cell $J$, the $\bbZ$-span of the $t_w$ with $w\in \H$ is closed under multiplication and forms an idempotent subalgebra $\sfA_\H$ of $\sfA_\J$. 

\begin{example}\label{b2asympalg}
Recall $\cS$ in type $B_2$ from Example \ref{cellex}\ref{Scellex} and consider the $H$-cell consisting of $\{s,sts\}$. Then the multiplication table in the Hecke algebra, when restricted to this $H$-cell and taken modulo bigger cells, is given by
\[
\begin{array}{c|cc}
&b_s&b_{sts}\\
\hline
b_s&(q+q^{-1})b_{s}&(q+q^{-1})b_{sts}\\
(q+q^{-1})b_{sts}&b_{sts}&(q+q^{-1})b_{s}\\
\end{array}
\]
and, since the $\mathbf{a}$-function on this cell is $1$, the multiplication table of the associated asymptotic algebra  $\sfA_\H$ is given by
\[
\begin{array}{c|cc}
&\mathsf{t}_s&\mathsf{t}_{sts}\\
\hline
\mathsf{t}_s&\mathsf{t}_{s}&\mathsf{t}_{sts}\\
\mathsf{t}_{sts}&\mathsf{t}_{sts}&\mathsf{t}_{s}\\
\end{array}.
\]
\end{example}

\subsection{Asymptotic bicategories}\label{abicatsec}

Fix an $H$-cell $\H$ and consider $\cS_\H$ as defined in Section \ref{Hcellredsec}. Following \cite[Subsection 3.2]{MMMTZ2}, we define
\begin{align*}
\X^0_\H&=\add\{ \theta_w \langle -j \rangle \,\mid\, w\in \H, j \geq 0\}\\
\Y^0_\H&=\add\{ \theta_w \langle -j \rangle \,\mid\, w\in \H, j > 0\}\\
\end{align*}
where the additive closures are taken in $\cS_\H(\varnothing,\varnothing)$. The quotient $\X^0_\H/\Y^0_\H$ is an additive category which carries a natural (non-strict) monoidal stucture. Thus, we can define the {\bf asymptotic bicategory} $\cA_\H$ as the bicategory on one object $\varnothing$ whose endormorphism category $\cA_\H(\varnothing,\varnothing)$ is given by $\X^0_\H/\Y^0_\H$. By construction, this categorifies $\sfA_\H$ and as explained in \cite[Subsection 3.2]{MMMTZ2}, it is fiat by \cite[Theorem 5.2]{EW3} and biequivalent to Lusztig's asymptotic monoidal category \cite[Subsection 18.19]{Lu2}. Moreover, $\cA_\H$ is a fusion category  \cite[Proposition 3.4]{MMMTZ2}, meaning it is also (locally) semisimple.
If $W$ is a finite Weyl group, $\cA_\H$ is well-understood and its simple birepresentations have been classified by Etingof, Ostrik and others. For more details, see \cite[Section 8]{MMMTZ2}.

\begin{example}\label{b2asympcat}
Again, consider $\cS$ in type $B_2$ from Example \ref{cellex}\ref{Scellex} and the $H$-cell consisting of $\{s,sts\}$. Then, setting  the multiplication table in $\cS$, when restricted to this $H$-cell and taken modulo bigger cells, is given by
\[
\begin{array}{c|cc}
&\theta_s&\theta_{sts}\\
\hline
\theta_s&\theta_{s}\oplus\theta_{s}\langle-2\rangle&\theta_{sts}\oplus\theta_{sts}\langle-2\rangle\\
\theta_{sts}&\theta_{sts}\oplus\theta_{sts}\langle-2\rangle&\theta_{s}\oplus\theta_{s}\langle-2\rangle\\
\end{array}
\]
and, denoting by $\rA_w$ the equivalence class of $\theta_w$ in $\cA_\H$, the multiplication table of the associated asymptotic bicategory  $\cA_\H$ is thus given by
\[
\begin{array}{c|cc}
&\rA_s&\rA_{sts}\\
\hline
\rA_s&\rA_{s}&\rA_{sts}\\
\rA_{sts}&\rA_{sts}&\rA_{s}.\\
\end{array}
\]
In this example, $\cA_\H$ is equivalent to the monoidal category of $\bbZ/2\bbZ$-graded vector spaces.
\end{example}

In order to classify simple $2$-representations of $\cS$, we would like to relate simple $2$-representations of $\cS_\H$ to simple birepresentations of $\cA_\H$.

\subsection{Relating $\cS_\H$ and $\cA_\H$.}

One of the crucial observations in \cite{MMMTZ2} is that, morally speaking, the asymptotic bicategory $\cA_\H$ is biequivalent to the endomorphism bicategory of the cell $2$-representation $\bfC_\H$ of $\cS_\H$. However, given that $\bfC_\H(\varnothing)$ is a graded category and $\cS_\H$ a (locally) graded $2$-category, the endomorphism bicategory $\cE nd_{\cS_\H}(\bfC)$ is naturally graded. Meanwhile the asymptotic bicategory $\cA_\H$ is ungraded. We thus define 
\begin{align*}
\X^t_\H&=\add\{ \theta_w \langle -j \rangle \,\mid\, w\in \H, j \geq t\}\\
\Y^t_\H&=\add\{ \theta_w \langle -j \rangle \,\mid\, w\in \H, j >t\}.\\
\end{align*}
and denote the quotient $ \X^t_\H/\Y^t_\H$ by $\A_\H^t$.
Then the monoidal structure on $\X^0_\H/\Y^0_\H$ induces a monoidal structure on
\begin{equation*}
\A^\bbZ:=\bigoplus_{t\in \bbZ} \A^t,
\end{equation*} 
with nontrivial components of the composition given by $\A_\H^s \times \A_\H^t \to\A_\H^{s+t}$. We denote by $\cA_\H^\bbZ$ the one-object bicategory on object $\varnothing$, whose endomorphism category is given by $\A^\bbZ$.

Recall that for a bicategory $\cC$, the bicategory $\cC^{ \mathsf{op}}$ is obtained by reversing the direction of all $1$-morphisms. Then the relationship between the asymptotic bicategory and endomorphisms of the cell $2$-representation can be stated as follows.

\begin{prop}\cite[Proposition 6.11]{MMMTZ2}
Let $\H$ be a diagonal $H$-cell in $\cS$ and let $\bfC_\H$ be the cell $2$-representation for $\H$ of the $2$-category $\cS_\H$. Then 
\[\cE nd_{\cS_\H}(\bfC_\H)\simeq\cA_\H^{\bbZ, \mathsf{op}}.\]
\end{prop}

We can thus view $\bfC(\varnothing)$ as having an action of $\cS_\H$ on the left and an action of $\cA_\H^\bbZ$ on the right, which mutually commute. Moreover, by the double centraliser theorem, the injective endomorphisms of $\bfC$ with respect to the right action of $\cA_\H^\bbZ$ are precisely given by $\add\H$.

Setting up a Morita context, one can prove the following.

\begin{prop}\label{SHAHsmod}\cite[Proposition 7.8]{MMMTZ2}
There is a biequivalence between the $2$-category of graded simple $2$-representations of $\cS_\H$ with apex $\H$ and the bicategory of graded simple birepresentations of $\cA_\H^\bbZ$.
\end{prop}

We observe that graded simple birepresentations of $\cA_\H^\bbZ$ are in bijection with simple birepresentations of $\cA_\H$. Namely,  for a simple birepresentation $\bfM$ of $\cA_\H$, let $\bfM^\bbZ$ be  the graded simple birepresentation of $\cA_\H$ defined as follows. Objects of $\bfM^\bbZ(\varnothing)$ are formal shifts $(X,t)$ where $X\in \bfM(\varnothing)$ and $t\in \bbZ$, and the action of a $1$-morphisms $\rF$ in $\A^t$ on $(X,s)$ is given by  $\bfM^\bbZ(\rF)(X,s) =(\bfM(\rF\langle t \rangle) ,s+t)$, noting that, by construction, $\rF\langle t \rangle\in \A_H^0 = \cA_\H(\varnothing,\varnothing)$. It is immediate that this assignment is injective and, using the fact that $\cA_\H$ is a fusion category, in particular semisimple, one can show that every graded simple birepresentation of $\cA_\H$ can be obtained in this way.

\subsection{Main Theorem}

Combining $H$-cell reduction with Proposition \ref{SHAHsmod}, we obtain the following main theorem.
\begin{thm} \label{mainthm}
There is a bijection between the set of graded simple $2$-representations of $\cS$ up to equivalence and the set of simple birepresentations of the $\cA_\H$ where $\H$ runs over a choice of diagonal $\H$-cell in every two-sided cell.
\end{thm}

\begin{rem} 
As remarked in Section \ref{abicatsec}, in most cases, simple birepresentations of the $\cA_\H$ have been classified. The only exceptions are Coxeter types $H_3$ and $H_4$, where for some two-sided cells, $\cA_\H$ is not understood for any diagonal $H$-cell.
Again, we refer to \cite[Subsection 8]{MMMTZ2} for more explanations and a complete list of cases.
\end{rem}

\begin{rem} 
The classification in Theorem \ref{mainthm} is constructive in principle, though in practice the computations involved will usually be difficult. Some explicit information that can be extracted in general is given in \cite[Section 5]{MMMTZ2}.
\end{rem}


{\bf Acknowledgements.}
I would like to thank my coauthors on the various papers contributing to this survey article for both their collaboration and their friendship.
I would also like to thank the referee for their helpful comments.

{\bf Funding.}
This work was partially supported by EPSRC grant EP/Z533750/1.


\end{document}